\documentclass[12pt]{amsart}
\usepackage{amsmath, amssymb, amsfonts, amsthm}
\usepackage{xspace, amscd, graphics}

\theoremstyle{plain}

\newtheorem{main}{Theorem}

\newtheorem{theorem}{Theorem}[section]
\newtheorem{proposition}[theorem]{Proposition}
\newtheorem{corollary}[theorem]{Corollary}
\newtheorem{lemma}[theorem]{Lemma}

\begin{document}

\title[Icosahedral Type]
{On the Symmetric Powers of Cusp Forms on $GL (2)$ of Icosahedral Type}
\date{Dec 30th, 2002}
\author{Song Wang}
\thanks{Partially supported by NSF grant \# 9729992}
\address{Department of Mathematics, \\
    Yale University, \\
    New Haven, CT 06520-8283.}
\email{song.wang@yale.edu}
\maketitle

\setcounter{section}{-1}
\section[Introduction]{Introduction}
\label{S:0}

In this Note, we prove three theorems.
Throughout, $F$ will denote a number field with
absolute Galois group $\mathcal{G}_{F} = Gal (\bar{F} / F)$,
and the adele ring $\mathbb{A}_{F} = F_{\infty} \times \mathbb{A}_{F, f}$.
When $\rho$ is an irreducible continuous $2$--dimensional
$\mathbb{C}$ representation of $\mathcal{G}_{F}$, one
says that it is {\it icosahedral}, resp.\ {\it octahedral},
resp.\ {\it tetrahedral}, resp.\ {\it dihedral} when
the projective image of $\rho (\mathcal{G}_{F})$ is $A_{5}$,
resp.\ $S_{4}$, resp.\ $A_{4}$, resp.\ $D_{2 m}$ for some $m \geq 1$.
Such $\rho$ is said to be {\it modular} if and only if there
exists a cuspidal automorphic representation
$\pi = \pi_{\infty} \otimes \pi_{f}$ of $GL_{2} (\mathbb{A}_{F})$
such that $L (s, \rho) = L (s, \pi_{f})$. Modularity is unknown
(in general) only when $\rho$ is icosahedral, in which case
$\rho$ is rational over $\mathbb{Q} (\sqrt{5})$.
Denoting by $\tau$ the nontrivial automorphism of $\mathbb{Q} (\sqrt{5})$,
we will say that
$\rho$ is {\it strongly modular} if both $\rho$ and $\rho^\tau$ are modular.
When
$F$ is totally real and $\rho$ totally odd, which is the primary case
of
interest, modularity implies strong modularity (see below). In the
Theorem
below, ${\rm sym}^{m} (\rho)$ denotes the symmetric $m$-th power of $\rho$,
i.e., the
composition of $\rho$ with the symmetric $m$-th power representation of
$GL_{2} (\mathbb{C})$ into $GL_{m + 1} (\mathbb{C})$.

\bigskip

The first main result is the following.

\medskip

\begin{main} \label{TM:A}

Let
\[
    \rho: \mathcal{G}_{F} \to GL_{2} (\mathbb{C})
\]
be a continuous irreducible, icosahedral representation
which is strongly modular, i.e., for which there exists a cuspidal
automorphic representation $\pi = \pi_{\infty} \otimes \pi_{f}$
of $GL_{2} (\mathbb{A}_{F})$, such that $L (s, \rho) = L (s, \pi_{f})$.
Then there exists a cuspidal representation
$\Pi = \Pi_{\infty} \otimes \Pi_{f}$
of $GL_{6} (\mathbb{A}_{F})$ such that
\[
    L (s, {\rm sym}^{5} (\rho)) = L (s, \Pi_{f})
\]
\end{main}

\bigskip

When $F = \mathbb{Q}$, many odd icosahedral
representations $\rho$ of $\mathcal{G}_{\mathbb{Q}} = {\rm Gal} (\bar{\mathbb{Q}} / \mathbb{Q}) $,
have been shown to be modular by R.~ Taylor, et al
(\cite{BDST}, \cite{Ta1}, \cite{Ta2}, \cite{BS});
The first example had been given by Buhler (\cite{Bu}).
In these cases the
associated
$\pi$ is generated by a holomorphic newform $\phi$ of weight $1$, and
$\phi^\tau$ is again a holomorphic newform of weight $1$. By a theorem
of
Deligne and Serre (\cite{D-S}), $\phi^\tau$ is associated to a
$2$-dimensional
representation $\rho_\tau$, which must be isomorphic to $\rho^\tau$ by
the
Chebotarev density theorem. Hence modularity implies strong
modularity for
odd icosahedral representations of $\mathcal{G}_{\mathbb{Q}}$. (The situation is
the
same when the base field $\mathbb{Q}$ is replaced by a totally real field $F$
as
long as $\rho$ is totally odd, and this is due to the analog of the
Deligne-Serre theorem due to Wiles (\cite{W}).)

\bigskip

By base change (\cite{AC}, \cite{La80}) we then get modular icosahedral
representations of $\mathcal{G}_{K}$ for any
cyclic extension $K$ of $\mathbb{Q}$.
So our theorem applies to these cases with no
hypothesis.

\bigskip

Given any $2$--dimensional irreducible icosahedral
representation $\rho$ of $\mathcal{G}_{F}$, one
sees that ${\rm sym}^{m} (\rho)$ is irreducible
if and only if $m \leq 5$ (see Section ~\ref{S:1}),
and the strong Artin conjecture, which is a part
of the Langlands philosophy, predicts the existence
of a cuspidal automorphic representation $\Pi_{m}$
of $GL_{m + 1} (\mathbb{A}_{F})$ for $m \leq 5$
with the same $L$--functions as ${\rm sym}^{m} (\rho)$.
When $\rho$ is strongly modular, the cuspidality of ${\rm sym}^{2} (\rho)$
has been known for a long time by the work
of Gelbart and Jacquet (\cite{GeJ79}),
and certain major recent works of H.~ Kim and F.~ Shahidi
(\cite{KSh2000}, \cite{KSh2001}, \cite{K2001})
establish this for $m = 3$ and $4$. In fact,
it is known by Kim (\cite{K2002}) when $F = \mathbb{Q}$ and $\rho$
odd
that every ${\rm sym}^{m} (\rho)$ is attached to an automorphic form on
$GL (m+1)$.
Briefly, for $m = 5$, the reason for this
is that ${\rm sym}^{5} (\rho)$ is twist equivalent to
$\rho \otimes {\rm sym}^{2} (\rho')$ where $\rho'$
is a Galois conjugate of $\rho$ (see Section ~\ref{S:1}).

\bigskip

Our main contribution here is to show that this
$\Pi$ here is indeed cuspidal on $GL (6) / F$. We prove it in two different ways.
One is to study the poles of the $L$--function,
and the other, which is perhaps of independent interest,
is to prove the following {\it cuspidality criterion}
for the Kim--Shahidi automorphic transfer from
$GL (2) \times GL (3)$ to $GL (6)$ (\cite{KSh2001}),
$(\pi, \eta) \mapsto \pi \boxtimes \eta$,
{\it when} $\eta$ is a twist of the symmetric square
of a cusp form on $GL (2)$.
More precisely, we prove the following:

\bigskip

\begin{main} \label{TM:B}
Let $\pi$, $\pi'$
be two cuspidal automorphic representations of
$GL_{2} (\mathbb{A}_{F})$, and let
$\Pi = \pi \boxtimes {\rm sym}^{2} (\pi')$ be the associated
isobaric automorphic representation of $GL_{6} (\mathbb{A}_{F})$.

Then $\Pi$ is cuspidal if and only if {\bf one} of the following
conditions hold:

\textnormal{(1)} $\pi = I_{K}^{F} (\chi)$ is dihedral for some
quadratic extension field $K$ of $F$,
$\pi'_{K}$, the base change of $\pi'$ to $K$, is not
dihedral, and
\[
    {\rm sym}^{2} (\pi'_{K}) \not\cong
    {\rm sym}^{2} (\pi'_{K}) \otimes \chi^{-1} (\chi \circ \theta)
\]
where $\theta$ is the nontrivial automorphism of $K / F$;

\textnormal{(2)} $\pi$ is not dihedral,
$\pi'$ is tetrahedral, or not of solvable polyhedral type,
and $Ad (\pi)$ and $Ad (\pi')$ are not equivalent.

\textnormal{(3)} $\pi$ is not dihedral,
$\pi'$ is octahedral, and $Ad (\pi)$ and $Ad (\pi')$
are not equivalent or twist equivalent by $\mu$
where $\mu$ is the global character corresponding
to the class field $K$ which is a quadratic extension field
of $F$ such that $\pi'_{K}$ is tetrahedral.

\end{main}

\bigskip

Recall that
$Ad (\pi) \cong {\rm sym}^{2} (\pi) \otimes \omega_{\pi}^{-1}$
where $\omega_{\pi}$ is the central character of $\pi$.
Also, note that if $\pi'$ is \emph{octahedral}, then $\mu$
is exactly the quadratic character such that
\[
    {\rm sym}^{3} (\pi') \cong {\rm sym}^{3} (\pi')
    \otimes \mu
\]

\medskip

We say that a cuspidal automorphic representation of
 $GL_{2} (\mathbb{A}_{F})$ is called not of \emph{solvable polyhedral type} (\cite{RaWa2001})
if it is not dihedral, tetrahedral or octahedral. It is a theorem of Kim--Shahidi that
if $\pi$ is not of solvable polyhedral type
then ${\rm sym}^{m} (\pi)$ is cuspidal for $m \leq 4$ (\cite{KSh2001}, \cite{K2001}).

\bigskip

Theorem ~\ref{TM:B} is proved in Section ~\ref{S:2}
below and the two proofs of Theorem ~\ref{TM:A}
will be given in Section ~\ref{S:3}.

\bigskip

Next, Recall that,
a Landau--Siegel zero of an $L$--function with a functional equation and Euler product
is a real zero of
this $L$--function close to $s = 1$
(see \cite{HRa95} and \cite{Ra99}).
Of course, the Generalized Riemann Hypothesis ({\it GRH})
implies the nonexistence of Landau--Siegel zeros for nice $L$--functions.
Unfortunately, this is obtained for only a few cases.
(\cite{HRa95}, \cite{RaWa2001}).

\medskip

Our third main result which will be proved in Section ~\ref{S:4}, is the following,
where we mean by a cusp form on $GL (2) / F$ of strongly icosahedral type
a cuspidal automorphic representation $\pi$
of $GL_{2} (\mathbb{A}_{F})$ attached
to a strongly modular icosahedral representation of $\mathcal{G}_{F}$:

\bigskip

\begin{main} \label{TM:C}
Let $\pi$ be a cusp form on $GL (2) / F$ of strongly icosahedral type, and
$\chi$ a idele character of $K$. Then
$L (s, {\rm sym}^{m} (\pi) \otimes \chi)$ has no Landau--Siegel
zero unless ${\rm sym}^{m} (\pi) \otimes \chi$
has a constituent of a trivial or quadratic
character $Q$. If this happens, $m$ is even, $Q =
\omega_{\pi}^{m/2} \chi^{m + 1}$ and there is at most one Landau--Siegel
zero which should come from the $L$--function of this character.
When $F = \mathbb{Q}$ and $\pi$ is self dual,
$L (s, {\rm sym}^{m} (\pi))$ has no Landau--Siegel zero at all.
\end{main}

\bigskip

\emph{Remark. } If $m < 12$, then the exceptional case will not happen
so that $L (s, {\rm sym}^{m} (\pi) \otimes \chi)$
has no Landau--Siegel zero. One can show by comparing central characters that
if ${\rm sym}^{m} (\pi)$ has a character as
its constituent, then $m$ is even
and this character
should be $\omega_{\pi}^{m/2}$.

\bigskip

This theorem needs more precise structure theory
(Theorem ~\ref{TM:D} in Section ~\ref{S:3})
of icosahedral representations
(see Section \ref{S:1} \& \ref{S:3},
 \cite{Bu}, \cite{FH}). The point (see also \cite{Bu}, \cite{K2002}) is that,
each twist of ${\rm sym}^{m} (\pi)$ is an isobaric sum of the twists
of the following (where $\pi^{\tau}$ is the Galois conjugate of
$\pi$ by $\tau$ which is the nontrivial element of
${\rm Gal} (\mathbb{Q} (\sqrt{5}) / \mathbb{Q})$):
\begin{align}
    &1, \,\pi, \,\pi^{\tau}, \,
    {\rm sym}^{2} (\pi), \,{\rm sym}^{2} (\pi^{\tau}),
    \notag \\
    &{\rm sym}^{3} (\pi), \,\pi \boxtimes \pi^{\tau}, \,{\rm sym}^{4} (\pi), \,
    {\rm sym}^{5} (\pi). \label{EQ:1} \tag{(A)}
\end{align}

So it suffices to show the nonexistence of Landau--Siegel zero
for the twist $L$--functions
of \eqref{EQ:1}. It is well known (\cite{HRa95}, \cite{Stk74})
that $L (s, \chi)$ has no Landau--Siegel zero unless $\chi$
is trivial or quadratic; The nonexistence of Landau--Siegel zero for
the twist $L$--functions for $\pi$ or $\pi'$
is obtained from \cite{HRa95};
From \cite{HRa95} and \cite{Ba97}, we also obtain the same for
the twist $L$--functions for ${\rm sym}^{2} (\pi)$
and ${\rm sym}^{2} (\pi^{\tau})$. From \cite{RaWa2001},
we get the same things for $\pi \times \pi^{\tau}$.
Furthermore, if $\pi$ is self dual or is twist equivalent to
a self dual automorphic representation, $L (s, {\rm sym}^{4} (\pi))$
has no Landau--Siegel zero (\cite{RaWa2001}).

\medskip

So we get
almost everything except for ${\rm sym}^{m} (\pi) \otimes \chi$
for $m = 3, 4$ or $5$. This is finally done by using a useful criterion
first formulated in \cite{HRa95} by D.~ Ramakrishnan
and J.~ Hoffstein (also developed in \cite{RaWa2001}),
the modularity for ${\rm sym}^{m} (\pi)$ (\cite{K2002}),
and Theorem ~\ref{TM:A}.

\bigskip

% Acknowledgement

This Note was inspired by a talk of
H.~ Kim based on \cite{K2002} at an MSRI conference
in Banff during 2001.
Of course, without the break-through works by Kim and Shahidi
on the functoriality on $GL (2) \times GL (3)$, ${\rm sym}^{3}$ (\cite{KSh2001})
and ${\rm sym}^{4}$ (\cite{K2001}), we cannot get these results.
Also, We would like to thank to D.~ Ramakrishnan and F.~ Shahidi
for the discussions during the preparation of this Note.
Finally, the author would like to thank IAS for membership during 2001--02,
as well as
NSF for the grant \# 9729992, and the department of mathematics at Caltech for providing a friendly place
and strong support during my visit in the summer of 2002.

\bigskip

\section[Structure Theory]
{Structure Theory for Icosahedral Representations}
\label{S:1}

In this section, we lay some facts
about icosahedral representations. Recall
that a Galois representation
$\rho: {\mathcal{G}}_{F} \to GL_{2} (\mathbb{C})$
is said to be \emph{icosahedral} if its image in
$PGL_{2} (\mathbb{C})$
is isomorphic to $A_{5}$.

\medskip

Let $\tilde{A}_{5}$ denotes the nontrivial central extension
of $A_{5}$ by $\mathbb{Z} / 2\mathbb{Z}$.
It is unique since
$H^{2} (A_{5}, \mathbb{Z} / 2\mathbb{Z}) \cong \mathbb{Z} / 2\mathbb{Z}$.
In fact, $\tilde{A}_{5} \cong SL_{2} (\mathbb{F}_{5})$.

\medskip

{

\begin{table}[h]
\caption{Character Table for $SL_{2} (\mathbb{F}_{5})$.}
 \label{TA:101}
\begin{tabular}{||c|c|c|c|c|c|c|c|c|c||}
    \hline
    \emph{\tiny Conj classes}
    & {\Tiny $\begin{pmatrix}  1 &  0 \\  0 &  1 \end{pmatrix}$}
    & {\Tiny $\begin{pmatrix} -1 &  0 \\  0 & -1 \end{pmatrix}$}
    & {\Tiny $\begin{pmatrix}  1 &  1 \\  0 &  1 \end{pmatrix}$}
    & {\Tiny $\begin{pmatrix}  1 &  2 \\  0 &  1 \end{pmatrix}$}
    & {\Tiny $\begin{pmatrix} -1 &  1 \\  0 & -1 \end{pmatrix}$}
    & {\Tiny $\begin{pmatrix} -1 &  2 \\  0 & -1 \end{pmatrix}$}
    & {\Tiny $\begin{pmatrix}  2 &  0 \\  0 &  3 \end{pmatrix}$}
    & {\Tiny $\begin{pmatrix}  3 &  2 \\  4 &  3 \end{pmatrix}$}
    & {\Tiny $\begin{pmatrix}  2 &  2 \\  4 &  2 \end{pmatrix}$}
    \\ \hline
    Size & $1$ & $1$ & $12$ & $12$ & $12$ & $12$ & $30$ & $20$ & $20$
    \\ \hline
    $U$ & $1$ & $1$ & $1$ & $1$ & $1$ & $1$ & $1$ & $1$ & $1$
    \\ \hline
    $V$ & $5$ & $5$ & $0$ & $0$ & $0$ & $0$ & $1$ & $-1$ & $-1$
    \\ \hline
    $W$ & $6$ & $6$ & $1$ & $1$ & $-1$ & $-1$ & $0$ & $0$ & $0$
    \\ \hline
    $X_{1}$ & $4$ & $-4$ & $-1$ & $-1$ & $1$ & $1$ & $0$ & $1$ & $-1$
    \\ \hline
    $X_{2}$ & $4$ & $4$ & $-1$ & $-1$ & $-1$ & $-1$ & $0$ & $1$ & $1$
    \\ \hline
    $W'$ & $3$ & $3$ & $\frac{1 + \sqrt{5}}{2}$ & $\frac{1 - \sqrt{5}}{2}$
     & $\frac{1 + \sqrt{5}}{2}$ & $\frac{1 - \sqrt{5}}{2}$ & $- 1$ & $0$ & $0$
    \\ \hline
    $W''$ & $3$ & $3$ & $\frac{1 - \sqrt{5}}{2}$ & $\frac{1 + \sqrt{5}}{2}$
     & $\frac{1 - \sqrt{5}}{2}$ & $\frac{1 + \sqrt{5}}{2}$ & $- 1$ & $0$ & $0$
    \\ \hline
    $X'$ & $2$ & $-2$ & $-\frac{1 + \sqrt{5}}{2}$ & $-\frac{1 - \sqrt{5}}{2}$
     & $\frac{1 + \sqrt{5}}{2}$ & $\frac{1 - \sqrt{5}}{2}$ & $0$ & $1$ & $-1$
    \\ \hline
    $X''$ & $2$ & $-2$ & $-\frac{1 - \sqrt{5}}{2}$ & $-\frac{1 + \sqrt{5}}{2}$
     & $\frac{1 - \sqrt{5}}{2}$ & $\frac{1 + \sqrt{5}}{2}$ & $0$ & $1$ & $-1$
    \\ \hline
\end{tabular}
\end{table}

}

From Table \ref{TA:101} (\cite{Bu}, \cite{NS}, \cite{FH}),
we see that there exist
two self dual irreducible $2$--dimensional representation of $\tilde{A}_{5}$, namely $X'$ and $X''$ which
are rational over $\mathbb{Q} (\sqrt{5})$.
Furthermore, $X'$ and $X''$
are conjugate by $\tau$ which is the nontrivial automorphism of
$\mathbb{Q} (\sqrt{5})/\mathbb{Q}$.
We use symbol $\rho_{\rm ico}$ for one of them, namely $X'$.
Hence we denote $\rho_{\rm ico}^{\tau}$ as $X''$.

\medskip

Also, the irreducible representations of $SL_{2} (\mathbb{F}_{5})$
are the following: (For
a proof the assertions, use the character table.)

\medskip

\begin{itemize}
    \item The trivial representation $U$;
    \item The $2$--dimensional representations $\rho_{\rm ico} = X'$,
    and $\rho_{\rm ico}^{\tau} = X''$,
      rational over $\mathbb{Q} (\sqrt{5})$;
    \item The $3$--dimensional representations
    $W' \cong {\rm sym}^{2} (\rho_{\rm ico})$,
      and $W'' \cong {\rm sym}^{2} (\rho_{\rm ico}^{\tau})$,
      which are rational over $\mathbb{Q} (\sqrt{5})$;
    \item The $4$--dimensional representation
    $X_{1} \cong {\rm sym}^{3} (\rho_{\rm ico})
      \cong {\rm sym}^{3} (\rho_{\rm ico}^{\tau})$,
    which is rational over $\mathbb{Q}$;
    \item The $4$--dimensional representation
    $X_{2} \cong \rho_{\rm ico} \otimes \rho_{\rm ico}^{\tau}$,
      also rational over $\mathbb{Q}$;
    \item The $5$--dimensional representation
        $V \cong {\rm sym}^{4} (\rho_{\rm ico})
         \cong {\rm sym}^{4} (\rho_{\rm ico}^{\tau})$,
        also rational over $\mathbb{Q}$;
    \item The $6$--dimensional representation
        $W \cong {\rm sym}^{2} (\rho_{\rm ico})
         \otimes \rho_{\rm ico}^{\tau}
         \cong {\rm sym}^{2} (\rho_{\rm ico}^{\tau})
         \otimes \rho_{\rm ico}$,
        also rational over $\mathbb{Q}$.
\end{itemize}

\medskip

More relations for the representations of $SL_{2} (\mathbb{F}_{5})$
(see also \cite{K2002}):

\begin{itemize}
    \item The symmetric $5$-the power of $\rho_{\rm ico}$,
        namely ${\rm sym}^{5} (\rho_{\rm ico})$,
        which is of dimension $6$, is equivalent
        to $W \cong {\rm sym}^{2} (\rho_{\rm ico})
        \otimes \rho_{\rm ico}^{\tau}
         \cong {\rm sym}^{2} (\rho_{\rm ico}^{\tau})
        \otimes \rho_{\rm ico}$;
        Also ${\rm sym}^{5} (\rho_{\rm ico}^{\tau})
        \cong {\rm sym}^{5} (\rho_{\rm ico}) \cong W$;
    \item The symmetric $6$-the power of $\rho_{\rm ico}$,
        namely ${\rm sym}^{6} (\rho_{\rm ico})$,
        is not irreducible, and is equivalent to
         ${\rm sym}^{2} (\rho_{\rm ico}^{\tau}) \oplus
         (\rho_{\rm ico} \otimes \rho_{\rm ico}^{\tau})$;
    \item The symmetric $7$-the power of $\rho_{\rm ico}$,
        namely ${\rm sym}^{7} (\rho_{\rm ico})$,
        is not irreducible either, and is equivalent to
         $\rho_{\rm ico}^{\tau} \oplus {\rm sym}^{5} (\rho_{\rm ico})$.
\end{itemize}

\bigskip

For the general icosahedral representation, we have the
following proposition:

\bigskip

\begin{proposition} \label{T:101}
Let $\rho$ be an icosahedral representation of $\mathcal{G}_{F}$,
and $G$ denotes $\rho (G)$, then $G$
is generated by its commutator subgroup
$G_{0}$ and its center $Z (G) \cong \mu_{2m}$,
which is a group of roots of unity of order $2 m$.
Furthermore, $G_{0}$ is isomorphic to $\tilde{A}_{5}$ with center
$\left\{\, \pm I \,\right\}$,
and $G \cong (G_{0} \times \mu_{2 m})/ \left\{\, \pm (I, 1) \,\right\}$.
Hence each irreducible
representation $\Lambda$ of $G$ can be expressed (uniquely)
as $(\Lambda_{0}, \mu)$ where $\Lambda_{0} = \Lambda |_{G}$
is an irreducible of $G_{0}$, and $\mu = \Lambda |_{\mu_{2m}}$
is a character of $\mu_{2m}$,
and such that $\Lambda_{0} (-I) = \mu (-1) I$.
Furthermore each such pair $(\Lambda_{0}, \mu)$ gives
an irreducible representation of $G$.
\end{proposition}

\bigskip

\emph{Remark: } Note that if $\rho$ is self dual of degree $2$,
$m = 1$, then $\rho$ is either the standard representation
or its Galois conjugation
by $\tau \in Aut (\mathbb{C})$ sending $\sqrt{5}$ to $- \sqrt{5}$.
Identify $G_{0}$ with $\tilde{A}_{5} \cong SL_{2} (\mathbb{F}_{5})$,
$\rho$ is $\rho_{\rm ico}$ or $\rho_{\rm ico}^{\tau}$.

\bigskip

\emph{Proof of Proposition ~\ref{T:101}.}

First consider the case
when ${\rm det} (G) = 1$,
i.e.\ ${\rm det} g = 1$ for all $g \in G$.
In this case, $G$ is a covering group of $A_{5}$ of degree $n$, where
$n = \# Z (G)$, while $Z (G) \subset Z (GL_{2} (\mathbb{C}))$.
As ${\rm det} (G) = 1$, $Z (G) \subset \{\, \pm I \,\}$.
Furthermore,
since $A_{5}$ has no irreducible representation of dimension $2$
(see \cite{FH}),
$Z (G)$ cannot be trivial. Thus $G$ is a nonsplit central
extension of $A_{5}$
by $\mathbb{Z} / 2 \mathbb{Z}$.
Thus $G \cong \tilde{A}_{5}$
(see the definition of $\tilde{A}_{5}$ at the beginning of this section).

\medskip

In general case, all elements of $G_{0} = (G, G)$,
which is the commutator group
of $G$, have determinant $1$,
and the image of $G_{0}$ in $PGL_{2} (\mathbb{C})$
is the same as the one of $G$,
and is also isomorphic to $A_{5}$
since $(A_{5}, A_{5}) = A_{5}$. We conclude that
$G = \left<\, G_{0}, Z (G) \,\right>
     \cong  (G_{0} \times \mu_{2m}) / \{\, \pm (I, 1) \,\}$
where $Z (G) \cong \mu_{2m}$ for some $m$.
From the discussion of the case ${\rm det} G = 1$,
we have $G_{0} \cong \tilde{A}_{5}$.
The proof of the rest assertions of the proposition
is then straight forward.

\qedsymbol

\bigskip

\begin{corollary} \label{T:102}
If $\Lambda_{1}$, $\Lambda_{2}$ are two representations of $G$
whose restrictions to $G_{0} \cong \tilde{A}_{5}$ are equivalent,
then they are
twist equivalent by a character. In fact, if $\Lambda_{i} = (\Lambda_{0}, \mu_{i})$,
then they are twist equivalent by $(1, \mu_{2} \mu_{1}^{-1})$ which is a character
of $G$ factoring through $\mu_{2m} / \{\, \pm 1\,\}$.
\end{corollary}

\qedsymbol

\bigskip

The following corollary describes all irreducible representations of $G$ and some relations.

\begin{corollary} \label{T:103}
Each irreducible representation of $G$ is twist equivalent
to one of the following:
\begin{align}
    1,\, &\Lambda_{\rm ico}, \,\Lambda'_{\rm ico},
    \, {\rm sym}^{2} (\Lambda_{\rm ico}),
    \, {\rm sym}^{2} (\Lambda'_{\rm ico}),
    \notag \\
    {\rm sym}^{3} (\Lambda_{\rm ico}),
    \, &{\rm sym}^{4} (\Lambda_{\rm ico}),
    \, {\rm sym}^{5} (\Lambda_{\rm ico}), \,
    \Lambda_{\rm ico} \otimes \Lambda'_{\rm ico}
    \notag
\end{align}
where $\Lambda_{\rm ico}$, $\Lambda'_{\rm ico}$ are two irreducible representations
of $G$ whose restrictions to $G_{0} \cong \tilde{A}_{5}$ are $\rho_{\rm ico}$
and $\rho_{\rm ico}^{\tau}$. Furthermore,
${\rm sym}^{m} (\Lambda_{\rm ico})$ and ${\rm sym}^{m} (\Lambda'_{\rm ico})$
are twist equivalent for $m = 3, 4$ and $5$;
${\rm sym}^{2} (\Lambda'_{\rm ico}) \otimes \Lambda_{\rm ico}$,
${\rm sym}^{2} (\Lambda_{\rm ico}) \otimes \Lambda'_{\rm ico}$,
${\rm sym}^{5} (\Lambda_{\rm ico})$ and ${\rm sym}^{5} (\Lambda'_{\rm ico})$
are twist equivalent.
\end{corollary}

\bigskip

\emph{Proof. } By Proposition ~\ref{T:101} and Corollary ~~\ref{T:102},
the first part is easy. For the
rest part, restrict all representations involved
to $G_{0}$, and apply Corollary ~\ref{T:102}.
We carry out the proof of the twist equivalence of ${\rm sym}^{5} (\Lambda_{\rm ico})$
and $\Lambda'_{\rm ico} \otimes {\rm sym}^{2} (\Lambda_{\rm ico})$
here, while the other relations are totally similar to deal with.

\medskip

Without loss of generality, say the restriction
of $\Lambda_{\rm ico}$ to $G_{0}$ is $\rho_{\rm ico}$.
Then the restrictions of
both sides are ${\rm sym}^{5} (\rho_{\rm ico})$
and $\rho^{\tau}_{\rm ico} \otimes {\rm sym}^{2} (\rho_{\rm ico})$
respectively.
They are equivalent from the discussion in this section.
Applying Corollary ~\ref{T:102}, we get the twist equivalence of
${\rm sym}^{5} (\Lambda_{\rm ico})$
and $\Lambda'_{\rm ico} \otimes {\rm sym}^{2} (\Lambda_{\rm ico})$.

\qedsymbol

\bigskip

Before concluding this section, we want to point out the following:

\bigskip

\begin{lemma} \label{T:104}
If $m$ is even, then $\Lambda_{\rm ico}$,
and $\Lambda'_{\rm ico}$ are not even twist equivalent
to a self dual representation;
Furthermore, there is no self dual $2$--dimensional
irreducible representation of $G$.
\end{lemma}

\bigskip

\emph{Proof. } each $2$--dimensional representation $\rho$ of $G$
is written as $(\rho_{0}, \mu)$.
Since $\rho_{0}$ is self dual, $\bar{\rho} = (\rho_{0}, \bar{\mu})$.
If $\rho$ is real, $\mu$  must be also real,
hence trivial or quadratic. However, as
$\mu (- 1) I = \rho_{0} (- I) = - I$
(since $\rho_{0}$ is either the standard representation
or its Galois conjugation by $\tau$),
and $\mu$ is of order divisible by $4$, hence $\mu$
cannot be real.

\qedsymbol

\bigskip

\section[Cuspidality Criterion]
{Cuspidality Criterion for $\pi \boxtimes {\rm sym}^{2} (\pi')$}
\label{S:2}

In this section, we will prove Theorem ~\ref{TM:B}.
Before this Note, no cuspidality criterion for the automorphic tensor product
(\cite{KSh2000}) on $GL (2) \times GL (3)$
was known.

\bigskip

Throughout this section, $\pi$ and $\pi'$ will be two cuspidal
automorphic representations of $GL_{2} (\mathbb{A}_{F})$.
First we will deal with the simpler case when $\pi$ is dihedral.

\medskip

\begin{lemma} \label{T:201}
Assume that $\pi = I_{K}^{F} (\chi)$ for some global character $\chi$
of $C_{K}$ where $K$ is a quadratic extension field of $F$.
Let $\Pi = \pi \boxtimes {\rm sym}^{2} (\pi')$.
Then
\[
    \Pi = I_{K}^{F} ({\rm sym}^{2} (\pi'_{K}) \otimes \chi)
\]
where $\pi'_{K}$ is the base change of $\pi'$ to $K$.

Hence $\Pi$ is cuspidal if and only if $\pi'_{K}$
is not dihedral, and
\[
    {\rm sym}^{2} (\pi'_{K}) \not\cong
    {\rm sym}^{2} (\pi'_{K}) \otimes \chi^{-1} (\chi \circ \theta)
\]
\end{lemma}

\emph{Remark: } As $\pi$ is cuspidal, $\chi \ne \chi \circ \theta$.
If $\pi'_{K}$ is not dihedral or tetrahedral, then of course
this lemma applies, as ${\rm sym}^{2} (\pi_{K})$ will not admit
a nontrivial self twist.

\bigskip

\emph{Proof. } The first statement is clear by the reciprocity
law of automorphic inductions and base changes (\cite{AC},
\cite{Cl}, \cite{HH}, \cite{La80}).
Let $\eta_{0} = {\rm sym}^{2} (\pi'_{K}) \otimes \chi$.
From Mackey's criterion
$\Pi = I_{K}^{F} (\eta_{0})$ is cuspidal
if and only if $\eta_{0}$ is cuspidal
and $\eta_{0} \not\cong \eta_{0} \otimes \chi^{-1} (\chi \circ \theta)$
where $\theta$ is the nontrivial automorphism of $K/F$.
Note that $\eta_{0}$ is cuspidal if and only if
$\pi'_{K}$ is not dihedral. Hence the second statement is true.

\qedsymbol

\bigskip

From now on, we assume that $\pi$ is not dihedral.
Of course, if $\pi'$ is dihedral, then ${\rm sym}^{2} (\pi')$
is not cuspidal and so is $\pi \boxtimes {\rm sym}^{2} (\pi')$.
Thus, Theorem ~\ref{TM:B} is finally reduced to the case
when $\pi$ and $\pi'$ are both nondihedral.

\bigskip

Let's recall a useful theorem about Rankin--Selberg $L$--functions
and isobaric decompositions
(\cite{JS90}, \cite{JPSS83}, \cite{La79-1}, \cite{La79-2}),

\bigskip

\begin{lemma} \label{T:202}
\textnormal{\bf (Jacquet--Shalika, Langlands)}

\textnormal{(1)} Let $\Pi$, $\tau$ be two automorphic representations
of $GL_{n} (\mathbb{A}_{F})$ and $GL_{m} (\mathbb{A}_{F})$
respectively. Assume that $\tau$ is cuspidal. Then
the order of the pole of $L (s, \Pi \otimes \tilde{\tau})$
is the same of the multiplicity of $\tau$ occurring in the isobaric
sum decomposition of $\Pi$.

\textnormal{(2)}
Let $\Pi$ be an isobaric automorphic representation of
$GL_{n} (\mathbb{A}_{F})$. Then $L (s, \Pi \times \tilde{\Pi})$
has a pole of order $m = \sum_{i} m_{i}^{2}$ at $s = 1$, where
$\Pi = \boxplus_{i} m_{i} \pi_{i}$ be the isobaric decomposition
of $\Pi$, and $\pi_{i}$ are inequivalent cuspidal representations
of smaller degree.

In particular, $m = 1$ if and only if $\Pi$ is cuspidal.
\end{lemma}

\qedsymbol

\bigskip

Now, we analyze $L (s, \Pi \otimes \tilde{\Pi})$ where
$\Pi = \pi \otimes {\rm sym}^{2} (\pi')$.
Let $\omega$ and $\omega'$ be the central characters of $\pi$
and $\pi'$.
Denote $Ad (\pi) = {\rm sym}^{2} (\pi) \otimes \omega^{-1}$,
$Ad (\pi') = {\rm sym}^{2} (\pi) \otimes {\omega'}^{-1}$
and $A^{4} (\pi') = {\rm sym}^{4} (\pi) \otimes {\omega'}^{-2}$.
Note that $Ad (\pi)$, $Ad (\pi')$ and $A^{4} (\pi')$ are self dual.

Hence, we have
\[
    \pi \boxtimes \tilde{\pi} = 1 \boxplus Ad (\pi)
\]
and
\begin{align}
    {\rm sym}^{2} (\pi') \boxtimes {\rm sym}^{2} (\tilde{\pi'})
    &= Ad (\pi') \boxtimes Ad (\pi')
    \notag \\
    &= 1 \boxplus Ad (\pi') \boxplus A^{4} (\pi')
    \notag
\end{align}

Thus
\begin{align}
    \Pi \boxtimes \tilde{\Pi}
    &= \pi \boxtimes \tilde{\pi}
    \boxtimes {\rm sym}^{2} (\pi') \boxtimes {\rm sym}^{2} (\tilde{\pi'})
    \notag \\
    &= (1 \boxplus Ad (\pi))
    \boxtimes (1 \boxplus Ad (\pi') \boxplus A^{4} (\pi'))
    \notag
\end{align}
Hence
\begin{align}
    L (s, \Pi \times \tilde{\Pi})
    &= \zeta_{F} (s) L (s, Ad (\pi)) L (s, Ad (\pi')) L (s, A^{4} (\pi'))
    \notag \\
    & \cdot L (s, Ad (\pi) \times Ad (\pi'))
    L (s, Ad (\pi) \times A^{4} (\pi'))
    \notag
\end{align}

Thus by Lemma ~\ref{T:202}, $\Pi$ is cuspidal if and only if
the order of the pole of $L (\Pi \times \tilde{\Pi})$
at $s = 1$ is $1$, if and only if the other $L$--factors
above other than $\zeta_{F} (s)$ are holomorphic at $s = 1$.
These lead to the following lemma: (Note that now $Ad (\pi)$
and $Ad (\pi')$ are cuspidal since $\pi$ and $\pi'$
are assumed to be nondihedral.)

\begin{lemma} \label{T:203}
If $\pi$ and $\pi'$ are not dihedral,
then $\Pi = \pi \boxtimes {\rm sym}^{2} (\pi')$ is cuspidal
if and only if {\bf all} the following hold:

\textnormal{(1)} $Ad (\pi)$ and $Ad (\pi')$ are not equivalent.

\textnormal{(2)} $A^{4} (\pi')$ does not have
     the trivial character as a constituent.

\textnormal{(3)} $A^{4} (\pi')$ does not have
     $Ad (\pi)$ as a constituent.
\end{lemma}

\qedsymbol

\bigskip

\begin{lemma} \label{T:204}
If $\pi'$ is not of solvable polyhedral type,
then (2) and (3) of Lemma ~\ref{T:203} hold.
\end{lemma}

\emph{Proof. }
From \cite{K2001}, ${\rm sym}^{4} (\pi')$ is cuspidal
and so is $A^{4} (\pi')$. Thus (2) and (3) of Lemma ~\ref{T:203}
hold.

\qedsymbol

\bigskip

\begin{lemma} \label{T:205}

\textnormal{(1)} If $\pi'$ is tetrahedral, then
\[
    A^{4} (\pi') \cong Ad (\pi') \boxplus \eta \boxplus \eta^{2}
\]
where $\eta$ is a cubic character such that
\[
    {\rm sym}^{2} (\pi') \cong {\rm sym}^{2} (\pi')
    \otimes \eta
\]

\textnormal{(2)} If $\pi'$ is octahedral, then
\[
    A^{4} (\pi') \cong Ad (\pi') \otimes \mu \boxplus \pi_{0}
\]
where $\mu$ is a quadratic character such that
\[
    {\rm sym}^{3} (\pi') \cong {\rm sym}^{3} (\pi')
    \otimes \mu
\]
and $\pi_{0} = I_{K}^{F} (\chi_{0})$ is some cuspidal
dihedral automorphic representation of $GL_{2} (\mathbb{A}_{F})$
where $K$ is the class field of $\mu$.
\end{lemma}

\emph{Proof. } See \cite{Tu81}, Theorem 3.3.7 of \cite{KSh2001}, and Section 5 of
\cite{RaWa2001}.

\qedsymbol

\bigskip

\emph{Proof of Theorem ~\ref{TM:B}.}

The case when $\pi$ is dihedral is dealt with in
Lemma ~\ref{T:201}. Now assume that $\pi$ and $\pi'$
are nondihedral.

\medskip

First we prove the necessity.
If $\Pi = \pi \boxtimes {\rm sym}^{2} (\pi')$ is cuspidal,
then Lemma ~\ref{T:203} and ~\ref{T:204} apply.
Hence $Ad (\pi)$ and $Ad (\pi')$ are not equivalent.
If $\pi'$ is octahedral with $\mu$ and $K$
described in Lemma ~\ref{T:205}, then $A^{4} (\pi')$
does not contain $Ad (\pi)$ as a constituent. Note that,
$Ad (\pi') \otimes \mu$ is a constituent of $A^{4} (\pi')$
hence it is not equivalent to $Ad (\pi)$. Then the necessity
is done.

\bigskip

Now the sufficiency. Assume first that $Ad (\pi)$
and $Ad (\pi')$ are not equivalent. (2)
and (3) of Lemma ~\ref{T:203} hold when $\pi'$ is not of
solvable polyhedral type. Then in this case, Lemma ~\ref{T:203}
applies, hence $\Pi$ is cuspidal.

\medskip

If $\pi'$ is tetrahedral, then from Lemma ~\ref{T:205}
the cuspidal constituents of $A^{4} (\pi')$
are $Ad (\pi')$ and two cubic characters. Hence (2)
of Lemma ~\ref{T:203} hold, and (1) and (3) are equivalent.
So the sufficiency in this case is proved.

\medskip

Finally, we deal with the case when $\pi'$
is octahedral. From Lemma ~\ref{T:205},
the only cuspidal constituent of $A^{4} (\pi')$
are $Ad (\pi') \otimes \mu$ and $I^{F}_{K} (\chi_{0})$.
So (2) of Lemma ~\ref{T:203} hold.
Thus if $Ad (\pi')$ and $Ad (\pi)$ are not equivalent
or twist equivalent by $\mu$, then
(1) and (3) hold, thus Lemma ~\ref{T:203}
applies. The sufficiency in this case is also obtained.

\medskip

Done.

\qedsymbol

\bigskip

\emph{Remark. } In fact, if $\pi'$ is octahedral,
$K$ is the quadratic field extension such that $\pi'_{K}$
is tetrahedral,
and $\chi_{0}$ is a cubic character of $C_{K}$
such that
\[
    {\rm sym}^{2} (\pi'_{K}) \cong {\rm sym}^{2} (\pi_{K})
         \otimes \chi_{0}
\]
then $A^{4} (\pi') = I^{F}_{K} (\chi_{0})$.

\bigskip

\section[Cuspidality of ${\rm sym}^{5} (\pi)$]
{Cuspidality of ${\rm sym}^{5} (\pi)$ for $\pi$ icosahedral}
\label{S:3}

In this section, we prove Theorem ~\ref{TM:A} in two
different ways.

\bigskip

Let $\rho$ be a strongly modular icosahedral representation
of $\mathcal{G}_{F}$ and $\pi$ the automorphic representation
of $GL_{2} (\mathbb{A}_{F})$ associated
with $\rho$. Then by the structure theory (Corollary ~\ref{T:103}),
${\rm sym}^{5} (\rho)$
is twist equivalent to $\rho^{\tau} \otimes {\rm sym}^{2} (\rho)$
and is irreducible.
The automorphy of ${\rm sym}^{5} (\rho)$
is known when $F = \mathbb{Q}$ and $\rho$ is odd (\cite{K2002}). One
immediately
gets the same for our $\rho$, and we indicate how. By assumption,
$\rho^\tau$ is modular, and since ${\rm sym}^{2} (\rho)$ is modular,
${\rm sym}^{2} (\rho) \otimes \rho^\tau$
and hence ${\rm sym}^{5} (\rho)$ is also modular by
\cite{K2002}.
\medskip

In view of this, Theorem ~\ref{TM:A} is a result of the following
known proposition which is an analogue of Lemma ~\ref{T:202}
on the Galois side.

\begin{proposition} \label{T:301}
If $\rho$ is an irreducible Galois representation
of $\mathcal{G}_{F}$, then $L (s, \rho \otimes \rho^{\vee})$
has a simple pole at $s = 1$. If $\rho$ is modular, and $\pi$
is the automorphic representation corresponding to
$\rho$, then $\pi$ is cuspidal if and only if $\rho$
is irreducible.
\end{proposition}

\emph{Proof. } (cf.\ Tate \cite{Tate})

One knows that given any $\mathbb{C}$--representation $\sigma$
of $\mathcal{G}_{F}$, we have
\[
    - {\rm ord}_{s = 1} L (s, \sigma)
    = {\rm dim}_{\mathbb{C}} {\rm Hom}_{\mathcal{G}_{F}} (1, \sigma^{\vee})
\]
Taking $\sigma$ to be $\rho \otimes \rho^{\vee}$, we see that the order of
pole is given by
\[
    {\rm dim}_{\mathbb{C}} {\rm Hom}_{\mathcal{G}_{F}} (1, \rho \otimes \rho^{\vee})
    = {\rm dim}_{\mathbb{C}} {\rm End}_{\mathcal{G}_{F}} (\rho)
\]
which, by Schur's lemma is $1$ if and only if $\rho$ is irreducible.

\medskip

Hence the first statement is clear.
In fact, for each Galois representation $\Lambda = \sum_{\tau} c_{\tau} \tau$
where $\tau$ are inequivalent irreducible representations of $\mathcal{G}_{F}$,
the order of pole of $L (s, \rho \otimes \rho^{\vee})$ at $s = 1$ is
$\sum_{\tau} C^{2}_{\tau}$.

\medskip

For the second part, we work with incomplete $L$--functions.
Let $S$ be a finite set of places of $F$
containing archimedean ones and the ones where $\rho$ (or $\pi$)
is ramified. Consider
\begin{align}
    L_{S} (s, \pi \times \tilde{\pi}) &=
    \prod_{v \in S} L (s, \pi_{v} \times \tilde{\pi}_{v})
    \notag \\
    L^{S} (s, \pi \times \tilde{\pi}) &=
    \prod_{v \notin S} L (s, \pi_{v} \times \tilde{\pi}_{v})
    \notag \\
    L_{S} (s, \rho \otimes \rho^{\vee}) &=
    \prod_{v \in S} L_{v} (s, \rho \otimes \rho^{\vee})
    \notag \\
    L^{S} (s, \rho \otimes \rho^{\vee}) &=
    \prod_{v \notin S} L_{v} (s, \rho \otimes \rho^{\vee})
    \notag
\end{align}
It is well known that each local $L$--factor $L (s, \pi_{v} \times \tilde{\pi}_{v})$
is holomorphic and not vanishing at $s = 1$ hence the order of the pole
of $L^{S} (s, \pi \times \tilde{\pi})$ is the same as $L (\pi \times \tilde{\pi})$
hence is $1$ if and only if $\pi$ is cuspidal from Lemma ~\ref{T:202}.

\medskip

Furthermore, for any Galois representation $\sigma$ and any nonarchimedean place $v$ of $F$,
$L_{v} (s, \sigma) = P (N \mathfrak{p}_{v}^{-s})^{-1}$ where $P$ is a polynomial with all
roots being of norm $1$. Hence $L_{v} (s, \sigma)$ is holomorphic and not vanishing
at $s = 1$. Thus the order of the pole of $L^{S} (s, \sigma)$ at $s = 1$ is exactly the
same as of $L (s, \sigma)$. Thus from the first statement of this proposition,
$L^{S} (s, \rho \otimes \rho^{\vee})$
has a simple pole if and only if $\rho$ is irreducible.

\medskip

Finally, as $\rho$ is modular, we have for all $v \notin S$,
\[
    L (s, \pi_{v} \times \tilde{\pi}_{v}) = L_{v} (\rho \otimes \rho^{\vee})
\]
Hence
\[
    L^{S} (s, \pi \times \tilde{\pi}) = L^{S} (\rho \otimes \rho^{\vee})
\]
Thus, $\pi$ is cuspidal if and only if
\[
    - {\rm ord}_{s = 1} L^{S} (s, \pi \times \tilde{\pi})
    = - {\rm ord}_{s = 1} L^{S} (s, \rho \otimes \rho^{\vee})
    = 1
\]
if and only if $\rho$ is irreducible.

\qedsymbol

\bigskip

\emph{Remark: } In general we don't know whether the following equality holds at {\bf ALL}
places $v$:
\[
    L (s, \pi_{v} \times \pi'_{v}) = L_{v} (s, \rho \otimes \rho')
\]
where $\rho$ and $\rho'$ are two modular Galois representations with
two automorphic representations $\pi$ and $\pi'$ associated to them respectively,
although we've known this
for those $v$ where $\rho_{v}$ and $\rho'_{v}$ are unramified.
When $\rho$ and $\rho'$ are $2$--dimensional, or one of $\rho$ and $\rho'$
is $2$--dimensional and the other one is $3$--dimensional, the automorphy
of $\rho \otimes \rho'$ (\cite{Ra2000}, \cite{KSh2000}) guarantees
this for {\bf all} $v$.

\bigskip

The second way to prove Theorem ~\ref{TM:A} is to apply
the criterion established in the previous section.
As we have seen, ${\rm sym}^{5} (\rho)$ is twist equivalent
to $\rho^{\tau} \otimes {\rm sym}^{2} (\rho)$, consequently,
${\rm sym}^{5} (\pi)$ is twist equivalent to
$\pi^{\tau} \boxtimes {\rm sym}^{2} (\pi)$,
where $\pi^{\tau} = \pi \circ \tau$ is the Galois conjugation of $\pi$
by $\tau$. Thus the condition (2)
of Theorem ~\ref{TM:B} holds hence this theorem applies.
In fact, $Ad (\pi^{\tau}) = {\rm sym}^{2} (\pi^{\tau}) \otimes
 \omega_{\pi^{\tau}}^{-1}$ and $Ad (\pi)
    = {\rm sym}^{2} (\pi) \otimes
 \omega_{\pi}^{-1}$ are not
equivalent as ${\rm sym}^{2} (\rho^{\tau})$
and ${\rm sym}^{2} (\rho)$ are not twist equivalent.

\bigskip

Now Theorem ~\ref{TM:A} is complete. Then we get a complete
structure theory for strongly modular icosahedral representations.

\bigskip

\emph{Notation: } Let $\pi$ be a cuspidal automorphic
representation of $GL_{2} (\mathbb{A}_{F})$ \emph{of strongly icosahedral type},
i.e., $\pi$ is associated to a strongly modular
icosahedral representation $\rho$ of $\mathcal{G}_{F}$.
Denote $M_{\rm ico} (\pi)$ as the set of
irreducible admissible representations generated by $\pi$ and $\pi^{\tau}$
via isobaric sums,
Rankin--Selberg products, twists and symmetric powers, where
$\pi^{\tau}$ is the cuspidal automorphic representation
of $GL_{2} (\mathbb{A}_{F})$ associated to $\rho^{\tau}$.

\medskip

\begin{main} \label{TM:D}

 \textnormal{(1)} all elements of $M_{\rm ico} (\pi)$
 are isobaric sums of the twists of the set $MG_{\rm ico} (\pi)$
 consisting of the following representations:
    \begin{align}
        &1, \,\pi, \,\pi^{\tau}, \,
        {\rm sym}^{2} (\pi), \,{\rm sym}^{2} (\pi^{\tau}),
        \notag \\
        &{\rm sym}^{3} (\pi), \,\pi \boxtimes \pi^{\tau},
     \,{\rm sym}^{4} (\pi), \,
        {\rm sym}^{5} (\pi).
        \notag
    \end{align}
    Furthermore, each two elements in $MG_{\rm ico} (\pi)$
    are not twist equivalent.
    ${\rm sym}^{m} (\pi)$ and ${\rm sym}^{m} (\pi^{\tau})$
    are twist equivalent for $m = 3$, $4$ and $5$. Also,
      ${\rm sym}^{5} (\pi)$, ${\rm sym}^{5} (\pi^{\tau})$,
    $\pi^{\tau} \boxtimes {\rm sym}^{2} (\pi)$ and
    $\pi \boxtimes {\rm sym}^{2} (\pi^{\tau})$ are twist equivalent.

 \textnormal{(2)} All elements in $MG_{\rm ico} (\pi)$ are automorphic. As a corollary,
    all elements in $M_{\rm ico} (\pi)$ are automorphic.

 \textnormal{(3)} All elements in $MG_{\rm ico} (\pi)$ are cuspidal.
\end{main}

\bigskip

\emph{Remark: } This theorem was first formulated by Kim
in \cite{K2002}.
The proof except for the cuspidality of ${\rm sym}^{5} (\pi)$
was also due to him.

\emph{Proof of Theorem ~\ref{TM:D}.}

Let $\rho = \rho_{\pi}$ be the odd icosahedral representation
associated
to $\pi$. Then $\rho$ and $\rho^{\tau}$ can be viewed
as representations of $G$ which is the image of $\rho$.
Hence $\rho = \Lambda_{\rm ico}$
or $\Lambda'_{\rm ico}$ (see Corollary ~\ref{T:103}), and all representations obtained from $\rho$
and $\rho^{\tau}$ via twists, direct sums,
tensor products and symmetric powers
are also viewed as representations of $G$. Then Corollary ~\ref{T:103}
applies, and thus (1) is proved.

\medskip

For (2) and (3), since $\pi$ and $\pi^{\tau}$ are
not of solvable polyhedral type,
${\rm sym}^{m} (\pi)$ and ${\rm sym}^{m}(\pi^{\tau})$
are cuspidal
for $m = 2$, $3$ and $4$ (\cite{GeJ79}, \cite{KSh2001}, \cite{K2001}).
Also, $\pi \boxtimes \pi^{\tau}$ is cuspidal (\cite{Ra2000}),
and $\pi^{\tau} \boxtimes {\rm sym}^{2} (\pi)$
is automorphic (\cite{KSh2001}) and cuspidal (Theorem ~\ref{TM:A}
or Theorem ~\ref{TM:B}).
Done.

\qedsymbol

\bigskip

Before we end this section, we would like to point out a result of
H.~ Kim in \cite{K2002} which asserts that ${\rm sym}^{4} (\pi)$
is monomial, thus giving an example of non-normal quintic
automorphic induction. Before this result, all known examples of
automorphic induction were for solvable extension (\cite{AC},
\cite{JPSS79}, \cite{Ha98} and \cite{Tu81}).

\bigskip

\begin{theorem} \label{T:302} \textnormal{\bf (H.~ Kim)}

Suppose that $K$ is an $A_{5}$--extension of $\mathbb{Q}$
satisfying the criteria as in \cite{BDST} or \cite{Ta2},
and that $\pi$ be a cuspidal automorphic
representation of strongly icosahedral type lifted from $K/\mathbb{Q}$.
Let $E$ be a non--normal quintic extension of
$\mathbb{Q}$ in $K$ such that ${\rm Gal} (K/E)$ is $A_{4}$. Let $N$ be
the unique cyclic cubic extension of $E$ in $K$. Let $\chi$ be the
global character of $C_{E}$ attached to the cubic extension $N/E$.

Then
${\rm I}^{\mathbb{Q}}_{E} (\chi)$ is equivalent to
$A^{4} (\pi) = {\rm sym}^{4} (\pi) \otimes \omega_{\pi}^{-2}$,
hence is automorphic.
\end{theorem}

\qedsymbol

\bigskip

\section[Landau--Siegel Zeros of $L (s, {\rm sym}^{m} (\pi) \otimes \chi)$]
{Landau--Siegel Zeros of $L (s, {\rm sym}^{m} (\pi) \otimes \chi)$}
\label{S:4}

In this section, the notations are the same as in the previous section.
Let us first quote the following useful criterion which is always
used for showing non--existence of Landau--Siegel zeros.

\begin{proposition} \label{T:401} \textnormal{(\cite{HRa95})}

Let $\pi$ be an isobaric automorphic representation
of $GL_{n} (\mathbb{A}_{F})$ with $L (s, \pi \times \bar{\pi})$ having
a pole of order $r \geq 1$ at $s = 1$. Then there is an effective
constant $c \geq 0$ depending on $n$ and $r$, such that $L (s, \pi \times \bar{\pi})$
has at most $r$ real zeros in the interval
\[
    J := \left\{\,
        s \in \mathbb{C} \,\mid\, 1 - c / \log M (\pi \times \bar{\pi}) < \Re (s)
        < 1
    \,\right\}.
\]
Furthermore, if $L (s, \pi \times \bar{\pi}) = {L_{1} (s)}^{k} L_{2} (s)$
for some nice $L$--series $L_{1} (s)$ and $L_{2} (s)$ with $k > r$,
and $L_{2} (s)$ holomorphic in $(t, 1)$ for some fixed $t \in (0, 1)$,
then $L_{1} (s)$ has no zeros in $J$.
\end{proposition}

\qedsymbol

\bigskip

\emph{Proof of Theorem ~\ref{TM:C}.}

From Theorem ~\ref{TM:D}, all ${\rm sym}^{m} (\pi) \otimes \chi$
are automorphic. Thus
it suffices to prove the nonexistence of Landau--Siegel zero
of $L (s, \Pi \otimes \chi)$, where $\Pi$ is $\pi$, $\pi^{\tau}$,
${\rm sym}^{2} (\pi)$, ${\rm sym}^{2} (\pi^{\tau})$,
${\rm sym}^{3} (\pi)$, ${\rm sym}^{4} (\pi)$,
${\rm sym}^{5} (\pi)$, or $\pi \boxtimes \pi^{\tau}$.

\medskip

If $\Pi \otimes \chi$ is not self dual, then it has no Landau--Siegel zero
(\cite{HRa95}). So
we need only to consider the case when $\Pi \otimes \chi$
is self dual.

\medskip

$L (s, \pi \otimes \chi)$, $L (s, \pi^{\tau} \otimes \chi)$
have no Landau--Siegel zero (\cite{HRa95}).

$L (s, {\rm sym}^{2} (\pi) \otimes \chi)$ has no Landau--Siegel zero. In
fact, when ${\rm sym}^{2} (\pi) \otimes \chi$ is self dual, its
central character is either trivial or quadratic. Thus, the
non-existence of Landau--Siegel zero follows from \cite{HRa95} and \cite{Ba97}.
When $\chi$ is trivial, we can
also get this from \cite{GHLL94}.

\medskip

$L (s, {\rm sym}^{m} (\pi) \otimes \chi)$
has no Landau--Siegel zero for $m = 3, 4, 5$.
This follows from the Lemma ~\ref{T:402}.

\medskip

(To be continued.)

\bigskip

\begin{lemma} \label{T:402}
Let $\pi$ be a nondihedral automorphic representation
of $GL (2)$ over $F$ such that ${\rm sym}^{m + 2} (\pi)$
and ${\rm sym}^{m - 2} (\pi)$
are automorphic and ${\rm sym}^{m} (\pi)$ are cuspidal
automorphic.
Then $L (s, {\rm sym}^{m} (\pi) \otimes \chi)$ has
no Landau--Siegel zero
for any Hecke character $\chi$ of $K$.
\end{lemma}

\bigskip

\emph{Proof of Lemma ~\ref{T:402}. }

Denote $\omega = \omega_{\pi}$ as the central character of $\pi$.

If ${\rm sym}^{m} (\pi) \otimes \chi$ is not self dual, then its $L$--function
has no Landau--Siegel zero.

\medskip

Now assume that ${\rm sym}^{m} (\pi) \otimes \chi$ is self dual.
Let $\Pi = 1 \boxplus ({\rm sym}^{m} (\pi) \otimes \chi)
     \boxplus ({\rm sym}^{2} (\pi) \otimes \omega^{-1})$,
then $\Pi$ is self dual, and $\Pi$ is an
isobaric sum of three cuspidal representations.
Hence $L (s, \Pi \times \Pi)$ has a pole of order $3$ at $1$.

\medskip

However,
\begin{align}
    L(s, & \Pi \times \Pi) =
    \zeta_{s}  (s) {L({\rm sym}^{m} (\pi) \otimes \chi)}^{2}
    {L (s, {\rm sym}^{2} (\pi) \otimes \omega^{-1})}^{2}
    \notag \\
    & \times {L (s, {\rm sym}^{2} (\pi) \otimes \omega^{-1}
        \times {\rm sym}^{m} (\pi) \otimes \chi)}^{2}
    L (s, {\rm sym}^{2} (\pi) \times {\rm sym}^{2} (\pi) \otimes \omega^{-2})
    \notag \\
    & \times L (s, {\rm sym}^{m} (\pi) \times {\rm sym}^{m} (\pi) \otimes \chi^{-2})
    \notag \\
    &= \zeta_{s} (s) {L({\rm sym}^{m} (\pi) \otimes \chi)}^{4}
    {L ({\rm sym}^{2} (\pi) \otimes \omega^{-1})}^{2}
    \notag \\
    & \times {L({\rm sym}^{m + 2} (\pi) \otimes \chi \omega^{-1})}^{2}
    {L({\rm sym}^{m - 2} (\pi) \otimes \chi \omega)}^{2}
    \notag \\
    & \times L (s, {\rm sym}^{2} (\pi) \times {\rm sym}^{2} (\pi) \otimes \omega^{-2})
    L (s, {\rm sym}^{m} (\pi) \times {\rm sym}^{m} (\pi) \otimes \chi^{-2})
    \notag
\end{align}
since
\begin{align}
    & {\rm sym}^{m} (\pi) \boxtimes {\rm sym}^{2} (\pi) =
    \notag \\
    & {\rm sym}^{m} (\pi) \boxplus {\rm sym}^{m + 2} (\pi) \otimes \omega^{-1}
    \boxplus {\rm sym}^{m - 2} (\pi) \otimes \omega
    \notag
\end{align}

Hence ${L (s, {\rm sym}^{m} (\pi) \otimes \chi)}^{4}$
divides
$L (s, \Pi \times \Pi)$, and the rest factors
are all automorphic $L$--functions.

\medskip

Thus, by Proposition ~\ref{T:401} (also \cite{HRa95}),
$L (s, {\rm sym}^{m} (\pi) \otimes \chi)$ has no Landau--Siegel zero.

\qedsymbol

\bigskip

\emph{Remark:} The non-existence of Landau--Siegel zero of
$L (s, {\rm sym}^{4} (\pi))$
when $\pi$ is self dual is followed from Theorem B of
\cite{RaWa2001}.
Unfortunately,
when $\pi = \pi (\rho)$ is a form corresponding
to an odd icosahedral representation,
$\pi$ cannot be self dual.

\medskip

\emph{Proof of Theorem ~\ref{TM:C} (Continued). }

Finally, the nonexistence of Landau--Siegel zero
for $\pi \boxtimes \pi^{\tau} \otimes \chi$
follows from Theorem A of \cite{RaWa2001} since $\pi$, $\pi^{\tau}$
are not dihedral and not twist equivalent. (Here
the form $\pi \boxtimes \pi^{\tau}$ is automorphic on $GL (4)$ (see \cite{Ra2000}).)

\medskip

The proof of the remaining statements are also straightforward.

\qedsymbol

\bigskip

\begin{proposition} \label{T:403}
Let $\pi$ be a cusp form on $GL (2)$ of strongly icosahedral type.
If a character $\omega'$ is an isobaric constituent of ${\rm sym}^{m} (\pi)$,
then $m$ is even and $\omega' = \omega_{\pi}^{m/2}$.
Also, ${\rm sym}^{m} (\pi)$ has no character
as its constituent for $m < 12$.
Hence $L (s, {\rm sym}^{m} (\pi) \otimes \chi)$
has no Landau--Siegel zero.
\end{proposition}

\bigskip

\emph{Proof. }

It is convenient to work on the Galois side. Let $\rho$ be an odd
icosahedral representation. Want to prove that if $\chi$ is
contained in ${\rm sym}^{m} (\rho)$, then $m$ is even,
and $\chi =
{\rm det} \rho^{m/2}$. In fact, writing $\rho = (\rho_{0}, \mu)$
as in Proposition ~\ref{T:101},
we have ${\rm sym}^{m} (\rho) = ({\rm
sym}^{m} (\rho_{0}), \mu^{m})$, and each irreducible component
should be $(1, \mu^{m})$ as $G_{0} \cong \tilde{A}_{5}$ has no
nontrivial $1$--dimensional representations (see Section
\ref{S:1}). Thus $m$ is even, and $(1, \mu^{m}) = {\rm det}
\rho^{m/2}$. Translating above to the automorphic side, and
noticing that $\omega_{\pi}$ is the global character corresponding
to ${\rm det} \rho$, we get the first statement.

\medskip

For the second, we need to verify the assertion for $m = 6$, $8$ and $10$.
We again work on the Galois representation side.
We want to prove that ${\rm sym}^{m} (\rho)$
has no constituent of character.
By the structure theory in Section \ref{S:1}, it suffices
to show that ${\rm sym}^{m} (\rho_{0})$
does not contain trivial representations.
This is true since
\[
    B = {({\rm sym}^{m/2} (\rho_{0}))}^{\otimes 2}
    = 1 \oplus (\oplus_{k = 1}^{m/2} {\rm sym}^{2 k} (\rho_{0})),
\]
which contains $1$
as multiplicity $1$ since ${\rm sym}^{m/2} (\rho_{0})$ is irreducible.
Hence ${\rm sym}^{m} (\rho_{0})$ cannot contain $1$.

\qedsymbol

\end{document}